\documentclass[letterpaper, 10 pt, conference]{ieeeconf}  

\IEEEoverridecommandlockouts                              

\usepackage{cite}
\usepackage{amsmath,amssymb,amsfonts}
\usepackage{algorithmic}
\usepackage{algorithm}
\usepackage{graphicx}
\usepackage{textcomp}
  \usepackage{subfigure,psfrag}

\usepackage{amssymb}
\usepackage{amsfonts, bm}

\newtheorem{rem}{Remark}
\newtheorem{defn}{Definition}
\newtheorem{lem}{Lemma}

\newtheorem{exmp}{Example}

  \usepackage{hyperref}

\title{\LARGE \bf
Data-driven approximation of control invariant set for linear system based on convex piecewise linear fitting
}

\author{Jun Xu$^{1}$ and Fanglin Chen$^{2}$
\thanks{*This work was not supported by any organization}
\thanks{$^{1}$Jun Xu is with School of Mechanical Engineering and Automation, Harbin Institute of Technology, Shenzhen, 518055, China, and also with the Key Laboratory of System Control and Information Processing, Ministry of Education, Shanghai 200240, China
        {\tt\small xujunqgy@hit.edu.cn}}%
\thanks{$^{2}$Fanglin Chen is with the Department of Computer Engineering, Harbin Institute of Technology, Shenzhen, 518055, China
        {\tt\small chenfanglin@hit.edu.cn }Corresponding author.}%
}

\begin{document}

\maketitle
\thispagestyle{empty}
\pagestyle{empty}

\begin{abstract}

Control invariant set is critical for guaranteeing safe control and  the problem of computing control invariant set for linear discrete-time system is revisited in this paper by using a data-driven approach. Specifically, sample points on convergent trajectories of linear MPC are recorded, of which the convex hull formulates a control invariant set for the linear system. To approximate the convex hull of multiple sample points, a convex piecewise linear (PWL) fitting framework has been proposed, which yields a polyhedral approximation with predefined complexity. A descent algorithm for the convex PWL fitting problem is also developed, which is guaranteed to converge to a local optimum. The proposed strategy is flexible in computing the control invariant set in high dimension with a predefined complexity. Simulation results show that the proposed data-driven approximation can compute the approximated control invariant set with high accuracy and relatively low computational cost.
\end{abstract}

\section{INTRODUCTION}

In scenarios for guaranteeing closed-loop stability of model predictive control (MPC) and safety-critical situations, the calculation of control invariant set is critical \cite{Mayne2000constrained,Tearle2021predictive}. A natural thought for guaranteeing safety or closed-loop stability is to drive the system states into the control invariant set, and the states are  guaranteed to be within the control invariant set since then. In recent years, this topic has gained considerable interests, especially in the area of safe learning for reinforcement learning \cite{Gros2020safe}.

Traditional method for computing the maximal control invariant set is based on geometric approach, in which the control invariant set is computed iteratively through the one-step backward operator until a termination condition is met \cite{Blanchini1999}.  The initial set can be either  the state constraint set or a control invariant set, corresponding to a non-increasing or non-decreasing sequence of sets, and the maximal control invariant set is approximated from outside or inside, respectively, \cite{Fiacchini2017computing,Rungger2017computing}.  In general, this kind of method is often prohibitive for high dimensional or complex systems, and it is difficult to guarantee the termination of the set-inclusion procedure \cite{Blanchini1999}.

An alternative procedure is to construct a more conservative but more computationally affordable sets, i.e., the control invariant set is assumed to be a 0-symmetric ellipsoid or polytope, and calculating the control invariant set then reduces to solve an optimization problem to maximize the volume of the ellipsoid or polytope. It has been proved in \cite{Pluymers2005efficient} that for linear systems with polyhedral constraints, the maximal control invariant set is also polyhedral, hence polyhedral approximation of the maximal control invariant set is more appropriate. Methods for calculating the polyhedral approximations fall into two categories. The first is by generating a set of non-decreasing control invariant set and introducing more flexible terminate conditions \cite{Gupta2021computation,Anevlavis2021controlled}. The second kind  is the data-based strategy, in which the safe set (actually control invariant) is constructed based on the closed-loop trajectory of the linear system. It was shown in \cite{Blanchini2003relatively} that the convex hull of states on a trajectory that terminates at the origin is control invariant. This idea is used in several work to construct safe sets or control invariant set, like  \cite{Brunner2013stabilizing,Rosolia2017learning,Wabersich2018linear}, etc. As a matter of fact, when the number of points on the trajectories is large, it is difficult to calculate the convex hull. In \cite{Wabersich2018linear}, a semi-definite programming problem was formulated to solve an approximated ellipsoid with maximal volume, which is computational expensive when the number of sample points is large.


In this paper, based on the idea of \cite{Blanchini2003relatively}, we propose a data-driven approach to approximate the maximal control invariant set. For linear systems with 0-symmetric constraints, we prove that 0-symmetric control invariant set exist. Then the convex hull of 0-symmetric convergent sample points of linear MPC are used to formulate  a 0-symmetric control invariant set. 
To approximate the convex hull of a large number of points in high dimension, we propose a convex piecewise linear (PWL) fitting framework, i.e., fitting the points with a convex PWL function while satisfying linear constraints.  An efficient procedure for removing redundant constraints is introduced to reduce the complexity of the convex PWL fitting problem and a descent algorithm  is put forward, which can converge to a local optimum of the optimization problem.

The rest of the paper is organized as follows. Section II describes some preliminaries, then sampling on the trajectories of linear MPC problem as well as the formulated control invariant set is introduced in Section III. Section IV describes the main work, i.e., polyhedral approximation of the control invariant set, in which the 0-symmetry of the control invariant set under symmetric constraints is proved, and the convex PWL fitting problem is proposed and solved. Simulation results are given in Section V and the paper ends with conclusions and future work in Section VI.

\section{Preliminaries}
Some basic definitions and assumptions are needed before proceeding to the computation of the control invariant set for linear systems. In this manuscript, we consider the constrained linear discrete-time system
\begin{equation}\label{eq:LTI}
\bm x_{k+1}=A\bm x_k+B\bm u_k, \bm x \in \mathcal{X}, \bm u \in \mathcal{U},
\end{equation}
in which $\bm x_k \in \mathbb{R}^{n_x}, \bm u_k \in \mathbb{R}^{n_u}$,  the set $\mathcal{X}$ and $\mathcal{U}$ are polyhedral and 0-symmetric. The definition 
of 0-symmetric convex set is given in Definition  \ref{def:0sym}.


\begin{defn}\label{def:0sym}
A  convex set $\mathcal{S}$ is 0-symmetric if
\[
\bm x \in \mathcal{S} \Rightarrow -\bm x \in \mathcal{S}.
\]
\end{defn}

The control invariant set for the constrained linear system (\ref{eq:LTI}) is defined as follows \cite{Kerrigan2000invariant}.
\begin{defn}
The non-empty set $\Omega \subset \mathbb{R}^n$ is a control invariant set for the system (\ref{eq:LTI}) if for all $\bm x_k \in \Omega$, there is $\bm u_k \in \mathcal{U}$, such that $\bm x_{k+1} \in \Omega$. The set $\Omega_{\max}$ is said to be maximal control invariant set for the system (\ref{eq:LTI}) if $\Omega \subset \Omega_{\max}, \forall \Omega$.
\end{defn}


\section{Sampling on the trajectories of linear MPC to formulate a control invariant set} \label{sec3}

According to \cite{Blanchini2003relatively}, the convex hull of points on convergent trajectories constitutes a control invariant set, and here we use linear MPC to generate such convergent trajectories.

\subsection{Linear MPC problem and the property of solutions}\label{sec:mpqp}
We first review the formulation of linear MPC problem and discuss the property of the solution.
The MPC problem for linear discrete-time time-invariant systems can be described as follows,
\begin{subequations}\label{eq:linearMPC}
\begin{align}
\min_U~&\bm x_{N_p}^TP\bm x_{N_p}+\sum_{k=0}^{N_p-1}(\bm x_k^TQ\bm x_k+\bm u_k^TR\bm u_k) \label{eq:obj_linearMPC}\\
s.t.~&\bm x_{k+1}=A\bm x_k+B\bm u_k, k=0, 1,\ldots,N_p-1, \label{eq:linear_dynamic}\\
{}&\bm x_k \in \mathcal{X}, k=1,\ldots,N_p, \label{eq:x_constraint}\\
{}&\bm u_k \in \mathcal{U}, k=0, 1,\ldots, N_p-1, \label{eq:u_constraint}\\
{}&\bm x_{N_p} \in \mathcal{X}_f,\label{eq:xn_constraint}
\end{align}
\end{subequations}
in which $\bm x_0\in \mathbb{R}^{n_x}$ is the initial state, and $N_p$ is the prediction horizon. The optimized variable is $U=[\bm u_0^T, \ldots, \bm u_{N_p-1}^T]^T\in \mathbb{R}^{n_u\cdot N_p}$ and $P, ~Q \in \mathbb{R}^{n_x \times n_x}, R \in \mathbb{R}^{n_u \times n_u}$ denote the weight matrices. The terminal set $\mathcal{X}_f$ is control invariant, i.e., once the state enters $\mathcal{X}_f$, there are feasible inputs such that the subsequent state remains in $\mathcal{X}_f$.

It was shown in \cite{Bemporad2002explicit} that the optimal solution $U^*$ is basically a continuous piecewise linear (PWL) function of $\bm x$, i.e., in subregions, $U^*$ admits an affine expression of $\bm x$, then the closed-loop system becomes a continuous PWL system, i.e.,
\begin{equation}\label{eq:closed_loop}
\bm x_{k+1}=A\bm x_k+Bf_{\rm CPWL}(\bm x_k),
\end{equation}
in which $f_{\rm CPWL}$ is a continuous PWL function.

It was shown in \cite{Mayne2000constrained} that if the terminal cost $\bm x_{N_p}^TP\bm x_{N_p}$ and the terminal set $\mathcal{X}_f$ are chosen suitably, the resulting closed-loop system (\ref{eq:closed_loop}) is exponentially stable within a region of attraction. It has also been proved in \cite{Scokaert1998constrained} that for stabilizing $(A, B)$, by increasing the prediction horizon $N_p$, one can obtain an optimal solution $U^*$ that can stabilize the linear system (\ref{eq:LTI}), by solving the linear MPC problem (\ref{eq:linearMPC}) with the constraint (\ref{eq:xn_constraint}) removed. Hence in this paper, we choose a suitable $N_p$ and omit the constraint (\ref{eq:xn_constraint}) to generate convergent trajectories, which is more computationally affordable as calculating the terminal invariant set requires additional work.

\subsection{Sample points collection}\label{sec:sampling}



In this subsection, {initial sample points are uniformly generated in $\mathcal{X}$}, according to Section \ref{sec:mpqp}, we can choose a long enough prediction horizon $N_p$ such that if the initial point is feasible, then the continuous PWL system (\ref{eq:closed_loop}) is stable and we record all the trajectories that converge to the equilibrium point $\bm 0$. 
Assume all the convergent points are,
\[
\mathcal{X}_s=\{\bm x_1,  \ldots, \bm x_{N_s}\},
\]
similar to \cite{Blanchini2003relatively},  we show that the convex hull of $\mathcal{X}_s$ is control invariant for the original linear system (\ref{eq:linear_dynamic}).

\begin{lem}\label{lem:control_invariant}
Suppose for the closed-loop system of linear MPC (\ref{eq:linearMPC}), every trajectory starting from $\bm x_i \in \mathcal{X}_s, i=1, \ldots, N_s$ converges to $\bm 0$, then the convex hull of $\bm x_1, \ldots, \bm x_{N_s}$, denoted as $\mathrm{conv}(\mathcal{X}_s)$ is control invariant for the linear system (\ref{eq:LTI}). Besides, if the vertices of the  polyhedral maximal control invariant set are all included in $\mathcal{X}_s$, the constructed convex hull $\mathrm{conv}(\mathcal{X}_s)$ is the maximal control invariant set for the linear system (\ref{eq:LTI}).
\end{lem}
\begin{proof}
For every $\bm x_i \in \mathcal{X}_s$, there is a control input $\bm u_i$, which is optimal for the MPC problem (\ref{eq:linearMPC}), such that
\[
A\bm x_i+B\bm u_i \in \mathrm{conv}(\mathcal{X}_s).
\]
Then for any $\bm x \in \mathrm{conv}(\mathcal{X}_s)$, i.e., there is a $\lambda$ such that
\begin{equation}\label{eq:cvx_comb}
\bm x=\sum_{i=1}^{N_s}\lambda_i \bm x_i,
\end{equation}
we can set $\bm u(\bm x)=\sum_{i=1}^{N_s}\lambda_i \bm u_i$, then we have
\[
A\bm x+B\bm u \in \mathrm{conv}(\mathcal{X}_s).
\]
Thus $\mathrm{conv}(\mathcal{X}_s)$ is control invariant.

If all the vertices of the polyhedral maximal control invariant set, say $\Omega_{\max}$, are included in $\mathcal{X}_s$, then from the definition of convex hull, we have
\[
\Omega_{\max} \subset \mathrm{conv}(\mathcal{X}_s).
\]
According to the definition of maximal control invariant set, we also have $\mathrm{conv}(\mathcal{X}_s) \subset \Omega_{
\max}$, which directly yields the conclusion that $\Omega_{\max}=\mathrm{conv}(\mathcal{X}_s)$.
\end{proof}
Then we demonstrate that a 0-symmetric controlled invariant set can always be obtained by 0-symmetric sample points.


\subsection{0-symmetric control invariant set}
\begin{lem}\label{lem:symmetric}
We can construct a  0-symmetric control invariant set $\Omega$, i.e., $\forall \bm x \in \Omega$, we have
\[
-\bm x \in \Omega,
\]
if the constraints set $\mathcal{X}$ and $\mathcal{U}$ are 0-symmetric. 
\end{lem}
\begin{proof}
According to the definition in \cite{Danielson2014symmetric}, if we choose two matrices $\Theta =-I_{n_x\times n_x}$ and $\Gamma =- I_{n_u \times n_u}$, in which $I_{n_x \times n_x}$ denotes the $n_x \times n_x$ identity matrix, then $(\Theta, \Gamma)$ is a symmetry of the dynamics (\ref{eq:linear_dynamic}), constraints (\ref{eq:x_constraint})-(\ref{eq:u_constraint}) and the cost (\ref{eq:obj_linearMPC}). Assume the set $\bar{\Omega}\subset \mathcal{X}$ is control invariant, i.e., for any $\bm x \in \bar{\Omega}$, there is some $\bm u \in \mathcal{U}$ such that $\bm x_{\rm next}=A\bm x+B\bm u \in \bar{\Omega}$.

According to \cite{Danielson2014symmetric}, as $\bm x$ is feasible, $-\bm x$ is also feasible. Denote $\tilde{\Omega}=\{-\bm x|\bm x \in \bar{\Omega}\}$, then we have
\[
-A\bm x-B\bm u=-\bm x_{\rm next} \in \tilde{\Omega}, \forall \bm x \in \bar{\Omega}.
\]
Let $\Omega = \bar{\Omega}\cup \tilde{\Omega}$, apparently $\Omega$ is 0-symmetric and for any $\bm x\in \Omega$ we have
\[
A\bm x+B\bm u \in \Omega,
\]
indicating that $\Omega$ is control invariant, and the conclusion follows.
\end{proof}

According to Lemma \ref{lem:control_invariant} and \ref{lem:symmetric}, if we have obtained the sample point set $\mathcal{X}_s$, we can augment the sample point set to obtain 0-symmetric sample point set, i.e., 
\begin{equation}\label{eq:sample_aug}
\tilde{\mathcal{X}}_s=\mathcal{X}_s \cup (-\mathcal{X}_s),
\end{equation}
To avoid the abuse of notation, hereafter we use $\mathcal{X}_s$ to represent $\tilde{\mathcal{X}}_s$. For any $\bm x\in \mathrm{conv}(\mathcal{X}_s)$, there is a $\lambda$ such that (\ref{eq:cvx_comb}) holds. Then we have
\[
-\bm x=\sum\limits_{i=1}^{N_s}\lambda_i (-\bm x_i),
\]
and $-\bm x  \in \mathcal{X}_s$, i.e., a 0-symmetric control invariant set can be obtained through the convex hull of ${\mathcal{X}}_s$.  







\section{Polyhedral approximation of the control invariant set}\label{sec:poly_approx}

Different from the existing research, which constructs a maximum-volume ellipsoid \cite{Wabersich2018linear}, here by investigating on the  0-symmetry of the control invariant set, we propose a method based on convex PWL fitting, i.e., fitting the convex hull of sample points, i.e., $\mathrm{conv}(\mathcal{X}_s)$, with a convex PWL function.


\subsection{Convex PWL fitting problem formulation}\label{sec:cvx_fitting}
As $\Omega$ is 0-symmetric, we only need to construct a polyhedral approximation of a part of $\Omega$. Specifically, we partition $\Omega$ into 3 sets, i.e., 
\begin{equation}\label{eq:partition}
\Omega=\Omega_0 \cup \Omega_{\rm N} \cup \Omega_{\rm P},
\end{equation}
in which the regions are defined as $\Omega_0=\Omega \cap \{\bm x|x_{n_x}=0\}$, and
\[
\begin{split}
\Omega_{\rm N} = \Omega \cap \{\bm x|x_{n_x}<0\},
\end{split}
\]
\[
\begin{split}
\Omega_{\rm P} = \Omega \cap \{\bm x|x_{n_x}>0 \},
\end{split}
\]
respectively. Here $x_j$ denotes the $j$-th component of $\bm x$, $j=1, \ldots, n_x$.

According to the partition (\ref{eq:partition}), all the sample points obtained in Section \ref{sec:sampling} are partitioned into 3 parts, i.e., we have index sets $\mathcal{I}_0, \mathcal{I}_{\rm N}$ and $\mathcal{I}_{\rm P}$, which are defined as,
\[
\mathcal{I}_0=\{i|\bm x_i \in \Omega_0\}, \mathcal{I}_{\rm N}=\{i|\bm x_i \in \Omega_{\rm N}\}, \mathcal{I}_{\rm P}=\{i|\bm x_i \in \Omega_{\rm P}\},
\]
and we have
\[
\{\bm x_i|i \in \mathcal{I}_0\} \cup \{\bm x_i|i \in \mathcal{I}_{\rm N}\} \cup \{\bm x_i|i \in \mathcal{I}_{\rm P}\}=\mathcal{X}_s.
\]

As is discussed before, the control invariant set $\Omega$ is 0-symmetric, a convex PWL function is constructed to formulate a lower bound for the points in $\Omega_0 \cup \Omega_{\rm N}$, the boundaries for points in $\Omega_0\cup \Omega_{\rm P}$ is then obtained by the property of symmetry. The detailed procedure is as follows.

 For the sample points $\bm x_i, i \in \mathcal{I}_0 \cup \mathcal{I}_{\rm N}$, 
a convex PWL fitting is constructed through designing a convex PWL function to approximate the lower bound of the convex hull $\mathrm{conv}(\bm x_i), i \in \mathcal{I}_0 \cup \mathcal{I}_{\rm N}$, i.e., 
\[\hat{x}_{n_x}=\bigvee\limits_{k \in \{1, \ldots, M\}} \ell_k(x_1, \ldots, x_{n_x-1}),
\] 
in which $\ell_k(x_1, \ldots, x_{n_x-1})$ is an affine function, such that 
\begin{equation}\label{eq:cvx_pwl}
x_{i,n_x}\geq \hat{x}_{i, n_x}, \forall i \in \mathcal{I}_0 \cup \mathcal{I}_{\rm N}.
\end{equation}
Here $x_{i, j}$ is the $j$-th component for sample point $\bm x_i$, $i=1,\ldots, N_s$, $j = 1, \ldots, n_x$, and $\hat{x}_{i, j}$ can be expressed as follows,
\begin{equation}\label{eq:hat_i}
\hat{x}_{i, n_x}=\bigvee_{k \in \{1, \ldots, M\}}\ell_k(x_{i,1}, \ldots, x_{i, n_x-1}).
\end{equation}


Then we have 
\begin{equation}\label{eq:inequality1}
\begin{split}
\ell_k(x_{i,1}, \ldots, x_{i, n_x-1})\leq x_{i,n_x}, \forall k \in \{1, \ldots, M\}, 
\forall i \in \mathcal{I}_0 \cup \mathcal{I}_{\rm N}.
\end{split}
\end{equation}

As for all $\bm x \in \Omega_0 \cup \Omega_{\rm P}$, the coordinate $x_{n_x}$ is larger than that for points $\bm x \in \Omega_0 \cup \Omega_{\rm N}$, hence (\ref{eq:inequality1}) is satisfied for all $\bm x_i \in \mathcal{X}_s$. Rearranging (\ref{eq:inequality1}) and we have the lower approximation to $\Omega$, i.e.,
\[
\mathcal{P}\bm x \leq q,
\]
where $\mathcal{P}$ and $q$ are obtained through (\ref{eq:inequality1}). 
According to the 0-symmetry of the control invariant set $\Omega$, we further have
\begin{equation}\label{eq:linear_inequality}
-\mathcal{P} \bm x \leq q,
\end{equation}
and the control invariant set $\Omega$ can be expressed as
\[
\left[\begin{array}{c}
\mathcal{P}\\
-\mathcal{P}\end{array}\right]\bm x\leq \left[\begin{array}{c}
q\\
q\end{array}\right].
\]

\begin{rem}
 The polyhedral approximation can also be constructed with respect to the points in $\Omega_0 \cup \Omega_{\rm P}$, i.e., $\hat{x}_{n_x}$ can also be constructed as an upper bound for the points with a  nonnegative $x_{n_x}$, i.e.,
\[
\hat{x}_{n_x}= \bigwedge \limits_{k \in \{1, \ldots, M\}} \ell_k(x_{1}, \ldots, x_{n_x-1}), \forall i \in \mathcal{I}_0 \cup \mathcal{I}_{\rm P},
\]
and the points $\bm x_i$ with $i \in \mathcal{I}_0 \cup \mathcal{I}_{\rm P}$ satisfies,
\[
x_{i,n_x}\leq \bigwedge \limits_{k \in \{1, \ldots, M\}} \ell_k(x_{i,1}, \ldots, x_{i,n_x-1}), \forall i \in \mathcal{I}_0 \cup \mathcal{I}_{\rm P}.
\]
This can also been used to get a linear inequality expression, which is actually (\ref{eq:linear_inequality}).
\end{rem}

The construction of the convex PWL fit to the sample points can be summarized as the following optimization problem, in which the distance between convex PWL function and the sample $x_{i, n_x}$ is minimized subject to the constraints that the sample points are above the convex PWL function. 

\begin{equation}\label{eq:cvx_fit_prob}
\begin{array}{rl}
\min\limits_{\ell_1, \ldots, \ell_M} &J=\sum_{i \in \mathcal{I}_0 \cup \mathcal{I}_{\rm N}}\left(\hat{x}_{i, n_x}-x_{i,n_x}\right)^2 \\
s.t. & (\ref{eq:hat_i})~ \mbox{and}~ (\ref{eq:inequality1}). 
\end{array}
\end{equation}
It is clear that for the above optimization problem,  the constraints are linear and the cost function is piecewise quadratic, i.e., 
\begin{equation}\label{eq:Jquadratic}
J=\sum\limits_{i \in \mathcal{I}_0 \cup \mathcal{I}_{\rm N}}\left(\ell_{k_i}(x_{i, 1}, \ldots, x_{i, n_x-1})-x_{i, n_x}\right)^2,
\end{equation}
for all $\ell_1, \ldots, \ell_M$ such that
 \begin{equation}\label{eq:linear_region}
 \begin{array}{r}
 \ell_{k_i}(x_{i,1}, \ldots, x_{i, n_x-1})\geq \ell_{k}(x_{i,1}, \ldots, x_{i, n_x-1}), \\
\quad\quad \forall i \in \mathcal{I}_0 \cup \mathcal{I}_{\rm N}, \forall k \in \{1, \ldots, M\}, k \neq k_i.
 \end{array}
 \end{equation}

\subsection{Reducing the complexity of the convex PWL fitting problem}
As can be seen from the optimization problem (\ref{eq:cvx_fit_prob}), the number of linear constraints is $M\cdot(|\mathcal{I}_{\rm N}|+|\mathcal{I}_0|)$, in which $|\cdot|$ denotes the number of elements in a set. When $|\mathcal{I}_{\rm N}|+|\mathcal{I}_0|$ is large, it is difficult to solve the optimization problem (\ref{eq:cvx_fit_prob}). Hence we perform a procedure to reduce the number of constraints and the following lemma provides the criterion for removing constraints. 

\begin{lem}\label{lem:remove}
If for some sample point $\bm x_{i_0}=[x_{i_0, 1}, \ldots, x_{i_0, n_x}]$, there exist sample points $\bm x_{j}, j\in \mathcal{J}_{i_0}$, where $\mathcal{J}_{i_0}$ is an index set with respect to $\bm x_{i_0}$, such that  $\bm x_{i_0}$ can be seen as the convex combination of $\bm x_{j}, j\in \mathcal{J}_{i_0}$, then if 
\[
\ell_k(x_{j,1}, \ldots, x_{j, n_x-1})\leq x_{j,n_x}, \forall k \in \{1, \ldots, M\}, \forall j \in \mathcal{J}_{i_0},
\]
we have
\[
\ell_k(x_{i_0,1}, \ldots, x_{i_0, n_x-1})\leq x_{i_0,n_x}, \forall k \in \{1, \ldots, M\}.
\]
\end{lem}

\begin{proof}
As  $\bm x_{i_0}$ is the convex combination of $\bm x_{j}, j\in \mathcal{J}_{i_0}$, there exist scalars $\lambda_j \in [0, 1], j \in \mathcal{J}_{i_0}$ such that
\[
\bm x_{i_0}=\sum\limits_{j \in \mathcal{J}_{i_0}} \lambda_j \bm x_j.
\]
Then $\forall k \in \{1, \ldots, M\}$, as $\ell_k$ is affine, we have
\[
\ell_k(x_{i_0,1}, \ldots, x_{i_0, n_x-1})=\sum\limits_{j \in \mathcal{J}_{i_0}} \lambda_j \ell_k(x_{j,1}, \ldots, x_{j, n_x-1}).
\]
Therefore, the following inequality holds,
\[
\begin{split}
&\ell_k(x_{i_0,1}, \ldots, x_{i_0, n_x-1})-x_{i_0,n_x}\\
&=\sum\limits_{j \in \mathcal{J}_{i_0}} \lambda_j \left(\ell_k(x_{j,1}, \ldots, x_{j, n_x-1})-x_{j,n_x}\right)\leq 0,
\end{split}
\]
and the conclusion follows.
\end{proof}

\begin{rem}
A simple way to eliminate redundant points is to formulate a simplex, which requires $n_x+1$ affinely independent points $P_1, \ldots, P_{n_x}, P_{n_x+1}$ in $n_x$-th dimension, here affinely independent means that $P_2-P_1, \ldots, P_{n_x+1}-P_1$ are linearly independent.
As the equilibrium point $\bm 0$ is in the centre of the control invariant set, it satisfies the constraint (\ref{eq:inequality1}). And if the distance between a samples points and $\bm 0$ is large, the sample point is likely to be on the boundary of the control invariant set. Hence we select $\bm 0$ and $n_x$ sample points with large distance to form a simplex, and the constraints corresponding to  points within the simplex can be removed.  The removal procedure can be repeated until the number of candidate points to form a simplex is less than $n_x+1$. 
\end{rem}


\subsection{Descent algorithm for solving the convex PWL fitting problem}
Assume $\tilde{\bm x}=[x_1, \ldots, x_{n_x-1}, 1]^T$, $\ell_k(x_1, \ldots, x_{n_x-1})=\bm \alpha_k^T \cdot \tilde{\bm x}$, in which \bm $\alpha_k \in \mathbb{R}^{n_x}$, then for each sample point $\bm x_i$, we have
\begin{equation}\label{eq:lk}
\ell_k(x_{i,1}, \ldots, x_{i, n_x-1})=\bm \alpha_k^T \cdot \tilde{\bm x}_i,
\end{equation}
in which $\tilde{\bm x}_i=[x_{i,1}, \ldots, x_{i, n_x-1}, 1]^T$.

For the convex PWL fitting problem (\ref{eq:cvx_fit_prob}), we propose to solve one quadratic programming at a time. 
Algorithm \ref{alg:pwq} shows the process of solving the convex PWL fitting problem, which is then shown descent and can converge to a local optimum in Lemma \ref{alg:pwq}.
  
 \begin{algorithm}
 \caption{Algorithm for solving convex PWL fitting problem (\ref{eq:cvx_fit_prob})}
 \label{alg:pwq}
 \begin{algorithmic}[1]
 \renewcommand{\algorithmicrequire}{\textbf{Input:} }
 \renewcommand{\algorithmicensure}{\textbf{Output:}}
\REQUIRE Optimization problem (\ref{eq:cvx_fit_prob}), sample points $\bm x_i, i \in \mathcal{I}_0 \cup \mathcal{I}_{\rm N}$, specified number of linear pieces $M$, a small enough threshold $\epsilon$.\\
\ENSURE  Solution  $\ell_1, \ldots, \ell_M$.
 \\
  \STATE Randomly generate $M$ vectors $\bm \alpha_{1,0}, \ldots, \bm \alpha_{M,0}$, let ${flag}=1$.
  \WHILE{$flag$}
  \STATE For each $i \in \mathcal{I}_0 \cup \mathcal{I}_{\rm N}$, find $k_{i,0}$ such that $\bm \alpha_{k_{i_0}}^T \cdot \tilde{\bm x}_i=\vee_{k \in \{1, \ldots, M\}}\bm \alpha_{k}^T \cdot \tilde{\bm x}_i$.
 \\ 
   \STATE Solve the following quadratic programming problem
   \begin{equation}\label{eq:qp}
   \begin{array}{rl}
   \min\limits_{\bm \alpha_1, \ldots, \bm \alpha_M} &\tilde{J} = \sum\limits_{i \in \mathcal{I}_0 \cup \mathcal{I}_{\rm N}} \left(\bm \alpha_{k_{i_0}}^T \cdot \tilde{\bm x}_i-x_{i,n_x}\right)^2 \\
   s.t. &  (\ref{eq:hat_i})~\mbox{and}~(\ref{eq:inequality1}) 
   \end{array}
   \end{equation}
   to obtain optimal $\bm \alpha_1, \ldots, \bm \alpha_M$. \label{line:qp}
     \IF{$J(\bm \alpha_1, \ldots, \bm \alpha_M)<J(\bm \alpha_{1,0}, \ldots, \bm \alpha_{M,0})$}
  \STATE Continue.
  \ELSE
  \STATE Using binary search in the line $[\mathcal{L}_0, \mathcal{L}]$, in which $\mathcal{L}_0=[\bm \alpha_{1,0}^T, \ldots, \bm \alpha_{M,0}^T]^T$ and $\mathcal{L}=[\bm \alpha_{1}^T, \ldots, \bm \alpha_{M}^T]^T$, to search for a point $\hat{\mathcal{L}}=[\hat{\bm \alpha}_{1}^T, \ldots, \hat{\bm \alpha}_{M}^T]^T$ such that
  \begin{equation}\label{eq:decrease_alg}
  J(\hat{\bm \alpha}_1, \ldots, \hat{\bm \alpha}_M)<J(\bm \alpha_{1,0}, \ldots, \bm \alpha_{M,0}).
  \end{equation}
     \vspace{-0.5cm}
  \STATE Let $\bm \alpha_1=\hat{\bm \alpha}_1, \ldots, \bm \alpha_M=\hat{\bm \alpha}_M$.
  \ENDIF
  \IF{$\|(\bm \alpha_1, \ldots, \bm \alpha_M)-(\bm \alpha_{1,0}, \ldots, \bm \alpha_{M,0})\|\leq \epsilon$}
  \STATE ${flag}=0$
  \ENDIF
  \STATE $\bm \alpha_{k, 0}=\bm \alpha_k, k=1, \ldots, M$.
  \ENDWHILE
 \RETURN $\ell_1, \ldots, \ell_M$ according to (\ref{eq:lk}). 
 \end{algorithmic} 
 \end{algorithm}

 \begin{lem}\label{lem:alg_convergent}
 Algorithm \ref{alg:pwq} is a descent algorithm and can converge to a local optimum of the optimization problem (\ref{eq:cvx_fit_prob}).
 \end{lem}
 \begin{proof}
 If Line 7 in Algorithm \ref{alg:pwq} is evaluated, a solution $(\bm \alpha_1, \ldots, \bm \alpha_M)$ with reduced value $J$ is obtained.
 
Then we demonstrate that if 
\[
J(\bm \alpha_1,\ldots, \bm \alpha_M)>J(\bm \alpha_{1,0}, \ldots, \bm \alpha_{M,0}),
\]
we can find some $\hat{\mathcal{L}}=\theta \mathcal{L}_0+(1-\theta) \mathcal{L}$ such that (\ref{eq:decrease_alg}) holds. This is easy to verify. For the line $[\mathcal{L}_0, \mathcal{L}]$, there is some $\hat{\mathcal{L}}$ such that 
\[
\hat{\bm \alpha}_{k_{i,0}}^T \cdot \tilde{\bm x}_i=\vee_{k \in \{1, \ldots, M\}}\hat{\bm \alpha}_k^T\cdot \tilde{\bm x}_i, \forall i \in \mathcal{I}_0 \cup \mathcal{I}_{\rm N},
\]
which means that for $\mathcal{L}_0$ and $\mathcal{L}$, $J$ admits the same quadratic function, i.e.,
\begin{subequations}\label{eq:proof_lem5}
\begin{align}
&J(\bm \alpha_{1,0},\ldots, \bm \alpha_{M,0})=\tilde{J}(\bm \alpha_{1,0},\ldots, \bm \alpha_{M,0}), \\
&J(\hat{\bm \alpha}_1, \ldots, \bm   \hat{\bm \alpha}_M)=\tilde{J}(\bm  \hat{\bm \alpha}_1, \ldots, \hat{\bm \alpha}_M).
\end{align}
\end{subequations}
According to the  optimization problem (\ref{eq:qp}), we have
\[
\tilde{J}(\bm \alpha_1, \ldots, \bm \alpha_M)<\tilde{J}(\bm \alpha_{1,0}, \ldots, \bm \alpha_{M,0}).
\]
As the function $\tilde{J}$ is quadratic, and $\hat{\mathcal{L}} \in [\mathcal{L}_0, \mathcal{L}]$, we further have
\[
\tilde{J}(\hat{\bm \alpha}_1, \ldots, \hat{\bm \alpha}_M)<\tilde{J}(\bm \alpha_{1,0}, \ldots, \bm \alpha_{M,0}).
\]
Combining with (\ref{eq:proof_lem5}), we have the following conclusion
\[
{J}(\hat{\bm \alpha}_1, \ldots, \hat{\bm \alpha}_M)<{J}({\bm \alpha}_{1,0}, \ldots, \bm {\alpha}_{M,0}),
\]
 meaning that Algorithm \ref{alg:pwq} is descent and can converge to a local optimum.
 \end{proof}

\begin{rem}
In our previous work \cite{Xu2012model},  a quadratic programming problem is also solved at a time, but within the linear subregion expressed by (\ref{eq:linear_region}). Compared with the approach proposed in this paper, additional $(M-1)\cdot (|\mathcal{I}_0|+|\mathcal{I}_{\rm N}|)$ linear inequalities are needed in \cite{Xu2012model} to describe the specific linear region when evaluating Line \ref{line:qp} in Algorithm \ref{alg:pwq}, resulting in a much larger computational complexity.
\end{rem}


In practice, we can set several candidate number of affine functions  $M$, which relates to the complexity of the approximated polyhedron. Then for each candidate $M$, several starting points are generated, resulting in a series of local optima. The final approximated polyhedron is the one chosen from all the local optima for all the candidate $M$ that yields the minimum $J$. It is noted that here the complexity of the polyhedral approximation can be predefined by tuning $M$.

\section{Simulation}\label{sec:simu}

In this section, we give two examples to illustrate the polyhedral approximation of the control invariant set, one is in two dimension and the other is in four dimension. 

\subsection{An example in two dimension}
\begin{exmp}\label{ex1}
Considering the two-dimensional example introduced in \cite{Danielson2012symmetric}, the system dynamics can be written as
\[
\begin{array}{rcl}
x_{k+1}&=&\left[\begin{array}{cc}
2&1\\
-1&2
\end{array}\right]x_k+\left[\begin{array}{cc}
1&0\\
0&1
\end{array}\right]u_k.
\end{array}
\]
 The  system constraints are $-1 \leq u_k \leq 1$ and $-1 \leq x_{k}\leq 1$. 
 
 
The parameters of the MPC problem used to generate sample points are  $Q=\mathrm{diag}(1,1)$, $R=5\times 10^3 \mathrm{diag}(1,1)$, $P=10$, and $N_p=10$. Uniformly generate 300 starting points in the domain $[-1, 1]^2$ and the sample points are generated on trajectories, which then yield 7044 samples in total, as Fig. \ref{fig:control_invariant} shows. Among all the 7044 sample points, only 33 are left after removing redundant constraints. Solving the convex PWL fitting problem for points with a nonnegative $x_2$ generates 5 affine functions, and the resulting polyhedral approximation contains 12 linear inequalities. The maximal control invariant set, which is computed through the MPT toolbox 3.0 \cite{MPT3}, the polyhedral approximation are shown in red and yellow, respectively. The processes of sampling, removing redundant constraints, and convex fitting cost 19.4477s, 0.97s, and 9.0012s, respectively. All the computations in this paper are implemented through MATLAB 2016b (MathWorks, USA) on an Apple M1 Max computer.
It can be seen from Fig. \ref{fig:control_invariant} that the approximated control invariant set approaches the genuine maximal control invariant set closely. 
 \begin{figure}[htbp]
\begin{center}
\includegraphics[width=0.7\columnwidth]{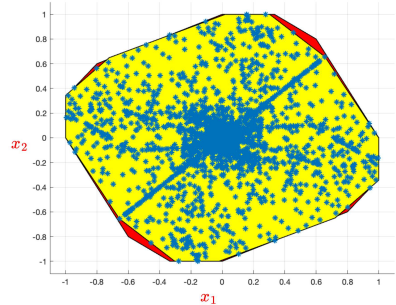}
\label{fig:control_invariant}
\caption{{Approximation of the control invariant set for Example 1. }}
\end{center}
\end{figure}



%
%
\end{exmp}

\subsection{A 4-dimensional example}
\begin{exmp}
Consider the following 4-dimensional inverted pendulum on cart example \cite{Karg2020efficient}, for which the system dynamic is,
\[
A=\left[
\begin{array}{cccc}
1&0.1&0&0\\
0&0.9818&0.2673&0\\
0&0&1&0.1\\
0&-0.0455&3.1182&1
\end{array}
\right], B=\left[\begin{array}{c}
0\\
0.1818\\
0\\
0.4546
\end{array}\right].
\]
The prediction horizon is taken to be $N=10$ and the value of the matrices in the cost function is $Q = \rm diag\{2, 2, 2, 2\}$, $R=1$, and $P = 1$. The system constraints are set to be $\|x\|_{\infty}\leq [1, 1.5, 0.35, 1]^T$ and $|u|\leq 1$. To get the convergent trajectories, we uniformly generate 300 starting points and record the system evolution, and then 0-symmetric points are generated as in (\ref{eq:sample_aug}) and we totally obtain 93954 sample points. After removing redundant points using Lemma \ref{lem:remove}, only 547 points are left for the convex PWL fitting problem (\ref{eq:cvx_fit_prob}), the result of which is a collection of 10 affine functions. Combined with the system constraints, we obtain an approximated control invariant set in 4-dimension with 26 linear inequalities constraints. Fig. \ref{fig:4d} shows the projection of the approximated control invariant set as well as the sample points. It can be seen from Fig. \ref{fig:4d} that the approximated control invariant set closely approximate the convex hull of  all the sample points.  
\begin{figure}[htbp]
\centering
\psfrag{x1}[c]{$x_1$}
\psfrag{x2}[c]{$x_2$}
\psfrag{x3}[c]{$x_3$}
\psfrag{x4}[c]{$x_4$}
\subfigure[Projection onto $x1-x2$ plane]{
\includegraphics[width=0.43\columnwidth]{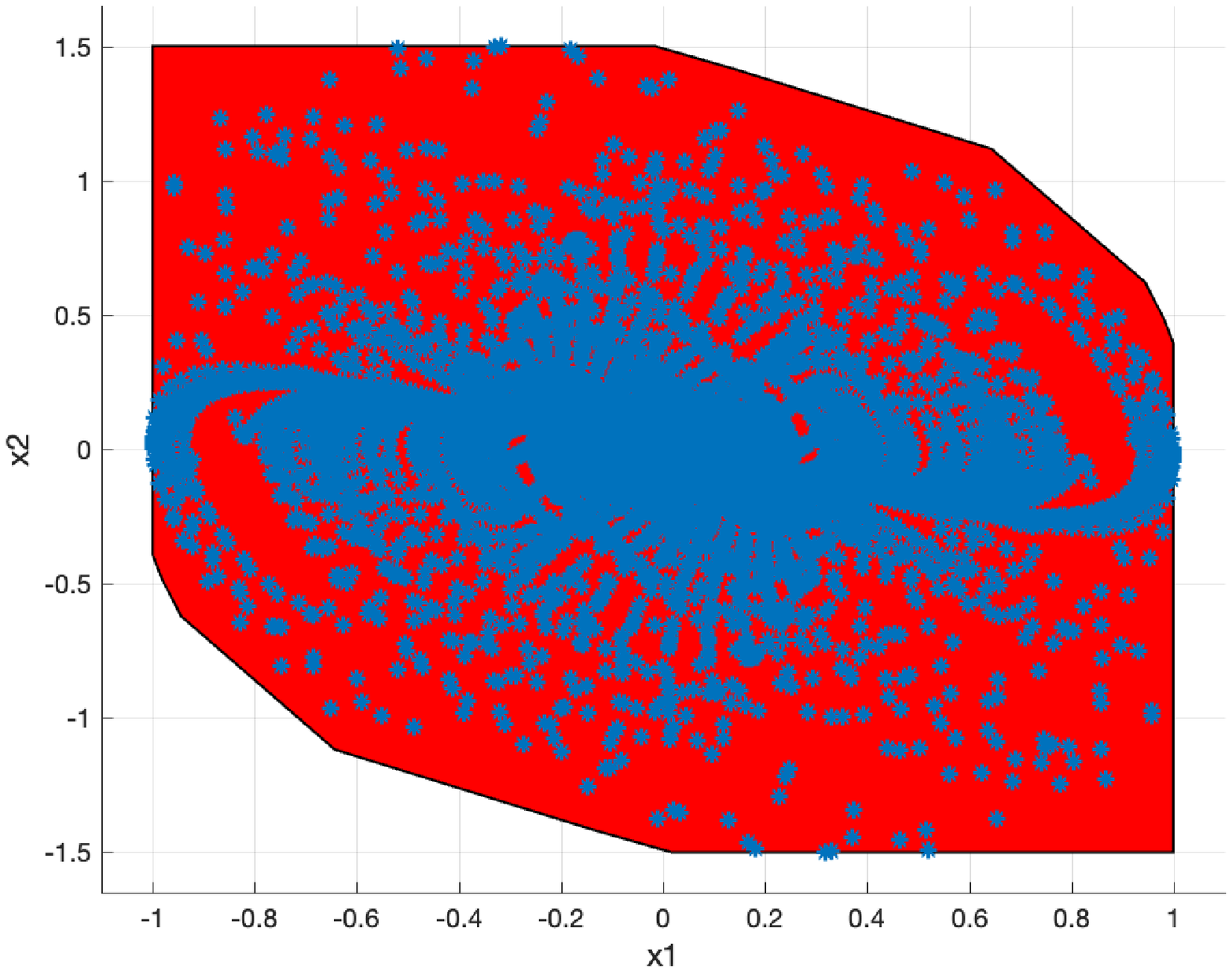}}
\subfigure[Projection onto $x1-x4$ plane]{
\includegraphics[width=0.43\columnwidth]{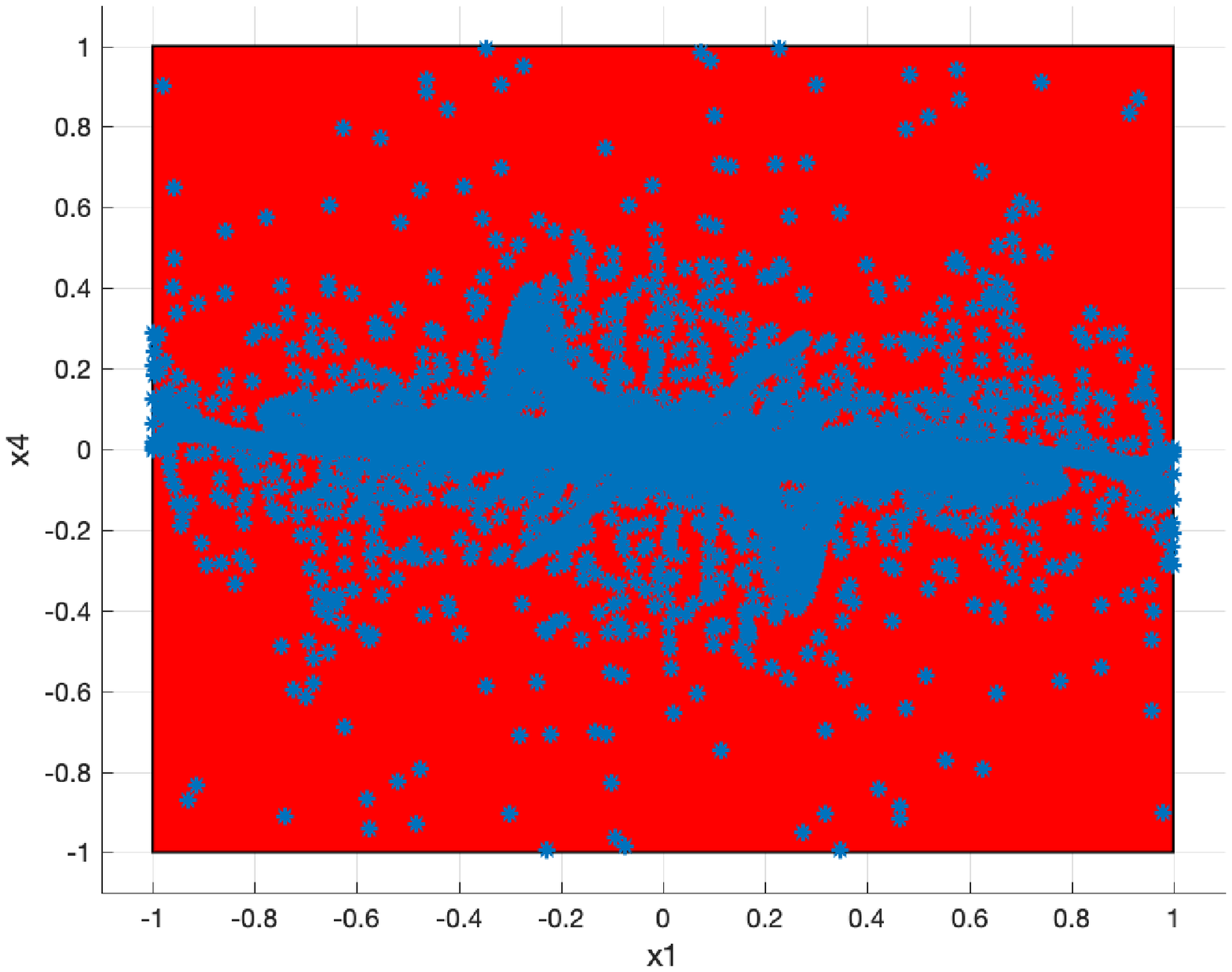}}
\subfigure[Projection onto $x2-x4$ plane]{
\includegraphics[width=0.43\columnwidth]{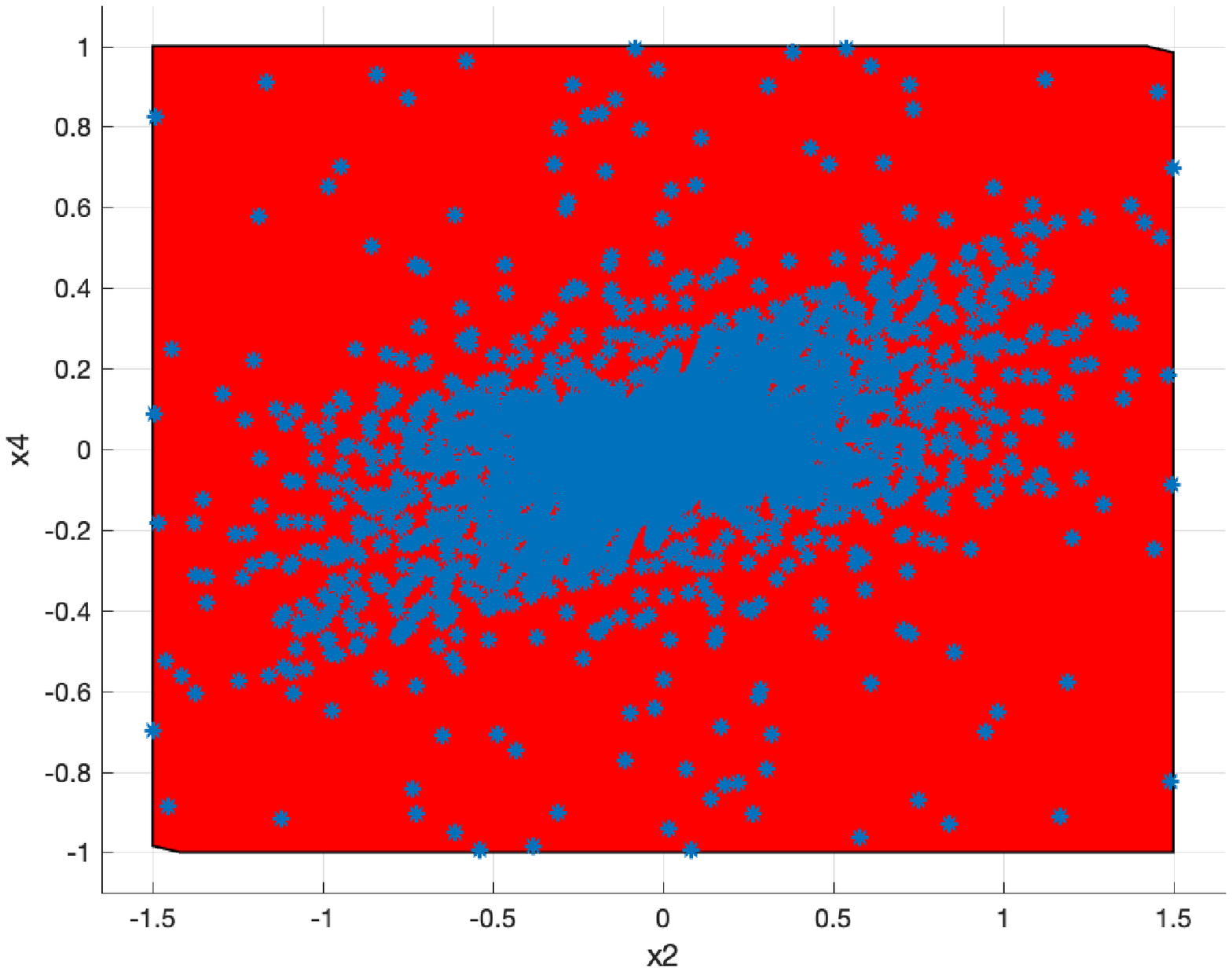}}
\subfigure[Projection onto $x3-x4$ plane]{
\includegraphics[width=0.43\columnwidth]{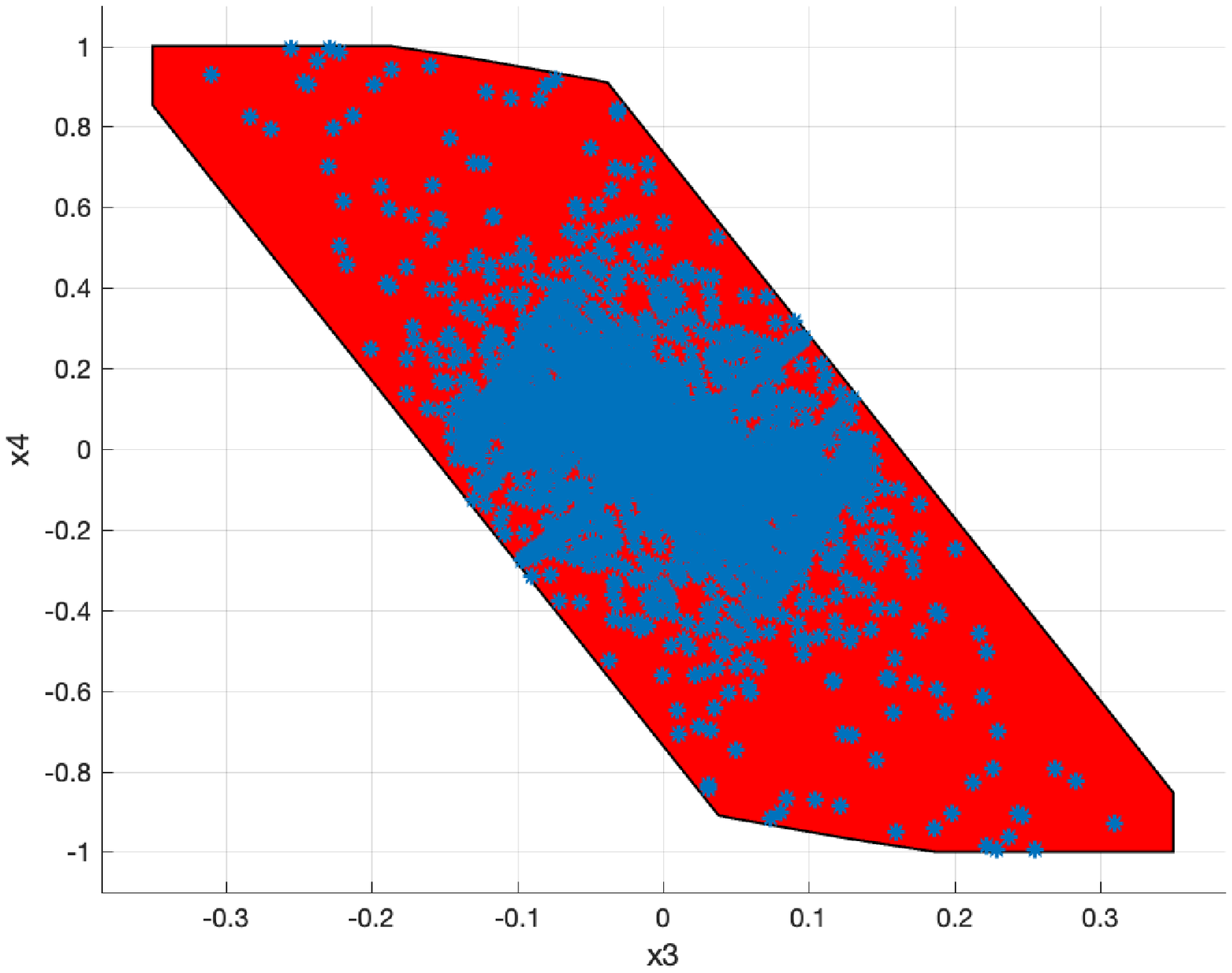}}
\caption{Projections of the approximated control invariant.}
\label{fig:4d}
\end{figure}

For this example, the sampling process and the constraints removing process takes 171.16s and 186.5s, respectively. Then setting the number of affine functions $M$ to be 10, 14,  18, and solving 100 convex PWL fitting problems with different starting points for each candidate $M$, the whole procedure costs 1070.1s, which is affordable for offline calculation. 
 It is noted that for this example, the MPT toolbox failed to calculate the exact control invariant set after 24 hours.
\end{exmp}

\section{CONCLUSIONS}

In this paper, we have developed a data-driven strategy for approximating the control invariant set of linear constrained system. If the state and input constraints are 0-symmetric, the control invariant set of a linear time-invariant system is proved to be 0-symmetric. Sample points on 0-symmetric convergent trajectories of linear MPC are collected, and a polyhedral approximation to the convex hull of all the sample points is constituted by formulating a convex PWL fitting problem. An algorithm for solving the convex PWL programming is provided, which is descent and can converge to the local minimum. Simulation results show the efficacy of the proposed strategy.

In the future, we will work towards constituting an approximated control invariant set that is  a subset of the maximal control invariant set.









\bibliographystyle{IEEEtran}
\bibliography{/Users/Jun/Documents/paper_work/refs}

\end{document}